\documentclass{article} 

\pdfoutput=1

\usepackage[latin1]{inputenc}
\usepackage[T1]{fontenc}
\usepackage[english]{babel}
\usepackage{times}
\usepackage{amssymb} 
\usepackage{graphicx} 
\usepackage{theorem}
\usepackage{float}

 \theoremstyle{plain} 
\theorembodyfont{\slshape} 
 \newtheorem{theorem}{Theorem}

{\theorembodyfont{\upshape}

\newtheorem{remark}[theorem]{Remark}}
\theoremheaderfont{\scshape\bfseries}

\newcommand{\titleref}[1]{\textsl{#1}}

\begin{document}

\title{\huge 
\textbf{A new LES model derived from generalized Navier-Stokes equations with nonlinear viscosity}}

\author{J.~M.~Rodr\'{\i}guez\thanks{Departamento de M\'etodos Matem\'aticos y de Representaci\'on, Universidade da Coru\~na, 15071 A Coru\~na (SPAIN). Email: jose.rodriguez.seijo@udc.es} ~and 
R.~Taboada-V\'azquez\thanks{Departamento de M\'etodos Matem\'aticos y de Representaci\'on, Universidade da Coru\~na, 15071 A Coru\~na (SPAIN). Email: raquel.taboada@udc.es}}

\date{}

\maketitle

\begin{abstract}
Large Eddy Simulation (LES) is a very useful tool when simulating turbulent flows if we are only interested in its ``larger'' scales (see \cite{Argyropoulos}, \cite{Guermond} or \cite{Sagaut}, for example). One of the possible ways to derive the LES equations is to apply a filter operator to the Navier-Stokes equations, obtaining a new equation governing the behavior of the filtered velocity. This approach introduces in the equations the 
so called \emph{subgrid-scale tensor}, that must be expressed in terms of the filtered velocity to close the problem. One of the most popular models is that proposed by Smagorinsky, where the subgrid-scale tensor is modeled by introducing an \emph{eddy viscosity}. 

In this work, we shall propose a new approximation to this problem by applying the filter, not to the Navier-Stokes equations, but to a generalized version of them with nonlinear viscosity. That is, we shall introduce a nonlinear viscosity, not as a procedure to close the subgrid-scale tensor, but as part of the model itself (see below). Consequently, we shall need a different method to close the subgrid-scale tensor, and we shall use the Clark approximation, where the Taylor expansion of the subgrid-scale tensor is computed (see \cite{Carati1} and \cite{Vreman}).

\vspace*{2ex}\par
\end{abstract}

\section{Introduction}

It is generally accepted that Navier-Stokes equations accurately model the behavior of incompressible viscous fluids. Simulating numerically these equations is possible today thanks to high computing power available. However, Direct Numerical Simulation (DNS) is limited to relatively low Reynolds numbers, because to simulate Navier-Stokes equations for a given value of Reynolds number ($R_e$) are needed, at least, $\mathcal{O}(R_e^{9/4})$ degrees of freedom (see \cite{Guermond}).

There are several methods to overcome this problem. For example, Reynolds Averaged Navier-Stokes models, $k-\varepsilon$ models or Large Eddy Simulation (LES). In this paper we shall talk about LES. Usually, LES is presented as an ``averaged'' or ``filtered'' version of Navier-Stokes equations.

Let us consider a ``filter'' operator $f \rightarrow \bar{f}$ (a spatial filter, a time filter or both), and let us assume that it is linear and commutes with spatial and time derivatives (see \cite{Sagaut} for different examples of this kind of filters). 
If we apply this filter to Navier-Stokes equations
\begin{eqnarray}
&&\frac{\partial \mathbf{u}}{\partial t} + \left ( \nabla \mathbf{u} \right ) \mathbf{u} 
+ \frac{1}{\rho_0} \nabla p - \nu \Delta \mathbf{u} = \mathbf{f}, \\
&&\nabla \cdot \mathbf{u} = 0,
\end{eqnarray}
(where $\mathbf{u}$ is the velocity field, $p$ is the pressure, and $\mathbf{f}$ is the density of external forces), we obtain 
\begin{eqnarray}
&&\frac{\partial \bar{\mathbf{u}}}{\partial t} + \left ( \nabla \bar{\mathbf{u}} \right ) \bar{\mathbf{u}} 
+ \frac{1}{\rho_0} \nabla \bar{p} - \nu \Delta \bar{\mathbf{u}} = \bar{\mathbf{f}} - \nabla \cdot \tau, \\
&&\nabla \cdot \bar{\mathbf{u}} = 0, \\
&&\tau_{ij} = \overline{u_i u_j} - \bar{u}_i\bar{u}_j, \label{def-tau}
\end{eqnarray}
where $\tau$ is the so called \emph{subgrid-scale tensor}. To close the LES model is necessary to express the subgrid-scale tensor only in terms of $\bar{\mathbf{u}}$. One of the most popular LES models is that proposed by Smagorinsky, where the subgrid-scale tensor is modeled by introducing an \emph{eddy viscosity} $\nu_e$ such that 
\begin{eqnarray} 
&&\tau = -2\nu_e\bar{\mathbf{D}}, 
\qquad \nu_e = (C_S \Delta x)^2 \sqrt{2}\, | \bar{\mathbf{D}} |, \label{eq-tau} \\
&&\bar{\mathbf{D}} = \frac{1}{2} \left ( \nabla \bar{\mathbf{u}} + \nabla \bar{\mathbf{u}}^T \right ), \qquad 
| \bar{\mathbf{D}} | = \left (\sum_{i,j=1}^3 \bar{D}_{ij} \bar{D}_{ij} \right )^{1/2}, \label{eq-tau-b}
\end{eqnarray}
$\Delta x$ is the subgrid-scale characteristic length 
and $C_S$ is a constant chosen to allow the model to emulate the kinetic-energy dissipation predicted by Kolmogorov (see \cite{Layton-libro}, \cite{Guermond} or \cite{Lesieur}).

The Smagorinsky model is sometimes too dissipative, and it has other theoretical problems, as the fact that subgrid-scale tensor in (\ref{eq-tau}) is odd in $\bar{\mathbf{u}}$, when by definition (see 
(\ref{def-tau})), it must be even in $\bar{\mathbf{u}}$. This happens because (\ref{eq-tau}) does not try to approximate the subgrid-scale tensor as defined in (\ref{def-tau}), but (as we have said before) to emulate the kinetic-energy dissipation predicted by Kolmogorov.

Interested reader can see other LES models in \cite{Layton-libro}, \cite{Sagaut}, \cite{Argyropoulos}, \cite{Lesieur} or \cite{Guermond}.

\section{Deriving the new LES model}

In the previous section, we have pointed out that the Smagorinsky model introduces dissipation through 
(\ref{eq-tau}), in a way that these equations are more a model than an approximation of the subgrid-scale tensor.

In this section, we shall propose to model this dissipation by introducing a nonlinear viscosity in the Navier-Stokes equations and then we shall apply the filter. Thus acting we introduce the dissipation predicted by Kolmogorov directly in our model and, when applying the filter, we shall close the LES model by approximating the subgrid-scale tensor using Clark approximation (see subsection \ref{subsection-filtering}).

\subsection{Introducing the nonlinear viscosity}

Let us consider the following generalization of Navier-Stokes equations
\begin{eqnarray}
&&\rho_0 \left ( \frac{\partial \mathbf{u}}{\partial t} + \left ( \nabla \mathbf{u} \right ) \mathbf{u} 
\right ) = \rho_0 \mathbf{f} + \nabla \cdot \mathbf{T}, \label{eq-NSNL-1}\\
&&\nabla \cdot \mathbf{u} = 0, \label{eq-NSNL-2}
\end{eqnarray}
where the stress tensor $\mathbf{T}$ is given by 
\begin{equation}
\mathbf{T} = -p\mathbf{I} + 2\mu_e\mathbf{D}, \qquad \mathbf{D} = \frac{1}{2} 
\left ( \nabla \mathbf{u} + \nabla \mathbf{u}^T \right ), \label{eq-T}
\end{equation}
and where the \emph{dynamic viscosity} $\mu_0$ has been substituted by the 
\emph{effective viscosity} $\mu_e$ depending on the norm of the strain rate tensor $\mathbf{D}$. 
To fix ideas, let us choose the following effective viscosity
\begin{equation}
\mu_e = \mu_e(| \mathbf{D} |) = \mu_0\left (1 + \lambda^2 | \mathbf{D} |^2 \right )^q \label{eq-mu-e}
\end{equation}
where $q>0$ and $\lambda > 0$. 

\begin{remark}
If we take $q=0$ or $\lambda = 0$, equations (\ref{eq-NSNL-1})-(\ref{eq-mu-e}) are the Navier-Stokes equations. Otherwise, they are known as Ladyzhenskaya model, and a unique global weak solution is guaranteed if $q \ge 1/4$ (see \cite{Guermond}). We obtain the Smagorinsky model if we take $q=1/2$, and properly choosing $\lambda$, we can recover the dissipation predicted by Kolmogorov.
\end{remark}

From the above remark, let us consider the effective viscosity given by (\ref{eq-mu-e})  with $q = 1/2$ in what follows. Thus, we introduce the dissipation predicted by Kolmogorov as part of the model, and we shall approach the subgrid-scale tensor in another way (see (\ref{eq-tau-21})-(\ref{eq-tau-22})).

\subsection{Filtering the equations} \label{subsection-filtering}

In this subsection, we apply a filter to the above equations ((\ref{eq-NSNL-1})-(\ref{eq-mu-e})). We are interested in linear filters that commute with spatial and time derivatives. Although there are several possible choices (see \cite{Sagaut}), we shall consider here a Gaussian filter:
\begin{equation}  
\bar{f}(t,\mathbf{x}) = \int_{-\infty}^{\infty} \int_{\mathbb{R}^3} G(s - t, \mathbf{y}- \mathbf{x})  f(s, \mathbf{y})\, \mathrm{d}\mathbf{y}\, \mathrm{d}s \label{eq-filtro-1}
\end{equation}
with 
\begin{equation}  
G(s, \mathbf{y}) = \frac{\gamma_T^{1/2}\gamma_L^{3/2}}{16\pi^2\eta^4} 
\exp\left ( - \frac{\gamma_T s^2 + \gamma_L |\mathbf{y}|^2}{4\eta^2} \right )
\label{eq-filtro-2}
\end{equation}
where $\eta>0$ is a small parameter related with the size of the filter, 
and $\gamma_T>0$, $\gamma_L>0$ are parameters related with the shape of the filter.

Following the ideas of Clark (see \cite{Carati1} and \cite{Vreman}), we approximate $f$ in (\ref{eq-filtro-1}) by its Taylor series, and then we can easily prove 
that this filter verifies for small values of $\eta$:

\begin{eqnarray}
&&\bar{f} = f + \eta^2 \left [ \frac{1}{\gamma_T}\frac{\partial^2 f}{\partial t^2} + 
\frac{1}{\gamma_L}\Delta f \right ]+ O(\eta^4), \label{eq-prom-1} \\
&&f = \bar{f} - \eta^2 \left [ \frac{1}{\gamma_T}\frac{\partial^2 \bar{f}}{\partial t^2} + 
\frac{1}{\gamma_L}\Delta \bar{f} \right ] + O(\eta^4), \\
&&\overline{fg} = \bar{f}\bar{g} + 2\eta^2\left [ \frac{1}{\gamma_T} 
\frac{\partial \bar{f}}{\partial t} \frac{\partial \bar{g}}{\partial t} + \frac{1}{\gamma_L} 
\nabla \bar{f} \cdot \nabla \bar{g} \right ] +  O(\eta^4), \\
&&\overline{fg} = \bar{f}g + \eta^2\Biggl [ \frac{2}{\gamma_T} \frac{\partial \bar{f}}{\partial t} 
\frac{\partial g}{\partial t} + \frac{1}{\gamma_L} 
\nabla \bar{f} \cdot \nabla g \\
&&\qquad \qquad {} + \bar{f}\left ( \frac{1}{\gamma_T}\frac{\partial^2 g}{\partial t^2} + 
\frac{1}{\gamma_L} \Delta g \right ) \Biggr ] + O(\eta^4), \nonumber 
\end{eqnarray}
and for any smooth enough real function $\beta$,
\begin{eqnarray}
&&\overline{\beta(f)g} = \beta(\bar{f})\bar{g} + \eta^2 \Biggl \{ 
2\beta'(\bar{f}) \left [ \frac{1}{\gamma_T} 
\frac{\partial \bar{f}}{\partial t} \frac{\partial \bar{g}}{\partial t} + \frac{1}{\gamma_L} 
\nabla \bar{f} \cdot \nabla \bar{g} \right ] \label{eq-prom-5} \\
&&\qquad {} + 
\beta''(\bar{f}) \left [ \frac{1}{\gamma_T} \left ( 
\frac{\partial \bar{f}}{\partial t} \right )^2 + \frac{1}{\gamma_L} 
\nabla \bar{f} \cdot \nabla \bar{f} \right ] \bar{g} \Biggr \} + O(\eta^4). \nonumber 
\end{eqnarray}

If we apply filter (\ref{eq-filtro-1})-(\ref{eq-filtro-2}) to equations (\ref{eq-NSNL-1})-(\ref{eq-mu-e}) (with $q=1/2$), 
and we take into account (\ref{eq-prom-1})-(\ref{eq-prom-5}), we obtain 
\begin{eqnarray}
&&\frac{\partial \bar{\mathbf{u}}}{\partial t} + \left ( \nabla \bar{\mathbf{u}} \right ) \bar{\mathbf{u}} 
+ \frac{1}{\rho_0} \nabla \bar{p} = \bar{\mathbf{f}} + \nabla \cdot \left (\mathbf{S} - \tau \right ), \\
&&\nabla \cdot \bar{\mathbf{u}} = 0, 
\end{eqnarray}
where $\tau$ is given by 
\begin{eqnarray}
&&\tau_{ij} = \overline{u_i u_j} - \bar{u}_i\bar{u}_j \label{eq-tau-21} \\
&&\quad = 2\left [ \frac{1}{\gamma_T} \frac{\partial \bar{u}_i}{\partial t} \frac{\partial \bar{u}_j}{\partial t} + 
\frac{1}{\gamma_L} \nabla \bar{u}_i \cdot \nabla \bar{u}_j \right ] \eta^2 + O(\eta^4), \label{eq-tau-22}
\end{eqnarray}
and $\mathbf{S}$ is given by 
\begin{eqnarray}
&&S_{ij} = \frac{2}{\rho_0}\left ( \overline{\mu_e(| \mathbf{D} |)\mathbf{D}} \right )_{ij} \label{eq-S} \\
&&\quad = 2\nu \Bigl [ \left (1 + \lambda^2 K_1 \right )^{1/2} 
+ \lambda^2 \left (1 + \lambda^2 K_1 \right )^{-1/2} K_2 \eta^2 \nonumber \\
&&\qquad {} - \frac{\lambda^4}{4} \left (1 + \lambda^2 K_1 \right )^{-3/2}K_3\eta^2 \Bigr ] \bar{D}_{ij} 
\nonumber \\
&&\quad {} + 2\lambda^2 \left (1 + \lambda^2 K_1 \right )^{-1/2} 
\hat{K}_{ij} \eta^2 + O(\eta^4) \nonumber
\end{eqnarray}
where we have introduced the notation 
\begin{eqnarray}
&&\nu = \mu_0/\rho_0, \quad K_1 = | \bar{\mathbf{D}} |^2, \label{eq-K-1} \\
&&K_2 = \sum_{m,n=1}^3 \left [ 
\frac{1}{\gamma_T} \left ( \frac{\partial \bar{D}_{mn}}{\partial t} \right )^2 + 
\frac{1}{\gamma_L} \nabla \bar{D}_{mn} \cdot \nabla \bar{D}_{mn}  \right ], \\
&&K_3 = \frac{1}{\gamma_T} \left ( \frac{\partial K_1}{\partial t} \right )^2 + 
\frac{1}{\gamma_L} \nabla K_1 \cdot \nabla K_1, \\
&&\hat{K}_{ij} = \frac{1}{\gamma_T} 
\frac{\partial K_1}{\partial t} \frac{\partial \bar{D}_{ij}}{\partial t} + 
\frac{1}{\gamma_L} \nabla K_1 \cdot \nabla \bar{D}_{ij}. \label{eq-K-4}
\end{eqnarray}

\begin{remark}
Formulas (\ref{eq-S})-(\ref{eq-K-4}) remember us those obtained in the dynamic procedure proposed by Germano (see \cite{Germano} and \cite{Lilly}), so this work could be seen as a new way to deduce a similar kind of models.
\end{remark}

\section{An example: application to Burgers equation}

In order to evaluate this model in the simplest possible case, we are going to apply, in this section, the same method to Burgers equation. The Burgers equation 
\begin{equation} \label{eq-Burgers}
\frac{\partial u}{\partial t} + u\frac{\partial u}{\partial x}  = \nu \frac{\partial^2 u}{\partial x^2} + f
\end{equation}
is often considered as an one-dimensional version of the Navier-Stokes equations. Next, 
we shall introduce a generalization of Burgers equation 
\begin{equation}
\rho_0\left ( \frac{\partial u}{\partial t} + u\frac{\partial u}{\partial x} \right ) = 
\rho_0 f \label{eq-Burgers-gen} 
+ \mu_0 \frac{\partial}{\partial x}\left [ \left (1 + \lambda^2 \left ( \frac{\partial u}{\partial x} 
\right )^2 \right )^{1/2} \frac{\partial u}{\partial x} \right ] 
\end{equation}
that can also be considered as an one-dimensional version of (\ref{eq-NSNL-1})-(\ref{eq-mu-e}).

\subsection{Filtering Burgers equation}

Let us now define the equivalent filter to (\ref{eq-filtro-1})-(\ref{eq-filtro-2}), but one dimensional in space, so we can apply it to the generalization of Burgers equation (\ref{eq-Burgers-gen}), 
\begin{equation}  
\bar{f}(t, x) = \int_{-\infty}^{\infty} \int_{-\infty}^{\infty} G(s - t, y-x)  f(s, y)\, \mathrm{d}y\, \mathrm{d}s \label{eq-filtro-1d-1}
\end{equation}
with 
\begin{equation}  \label{filtro-eq-Burgers} 
G(s, y) = \frac{\gamma_T^{1/2}\gamma_L^{1/2}}{4\pi\eta^2} 
\exp\left ( - \frac{\gamma_T s^2 + \gamma_L y^2}{4\eta^2} \right )
\label{eq-filtro-1d-2}
\end{equation}

If we now apply filter (\ref{eq-filtro-1d-1})-(\ref{eq-filtro-1d-2}) to equation (\ref{eq-Burgers-gen}), 
and we take into account the formulas corresponding to (\ref{eq-prom-1})-(\ref{eq-prom-5}) in the 
one-dimensional spatial case, we obtain

\begin{equation} \label{eq-Burgers-prom}
\frac{\partial \bar{u}}{\partial t} + \bar{u}\frac{\partial \bar{u}}{\partial x} = 
\bar{f} + \frac{\partial}{\partial x}\left ( S - \tau \right )
\end{equation}
where 
\begin{eqnarray}
&&\tau = \frac{1}{2}\left [ \overline{u^2} - (\bar{u})^2 \right ] 
= \left [ \frac{1}{\gamma_T} \left ( \frac{\partial \bar{u}}{\partial t} \right )^2 + 
\frac{1}{\gamma_L} \left ( \frac{\partial \bar{u}}{\partial x} \right )^2 \right ] \eta^2 + O(\eta^4), \\
&& \nonumber \\
&& \nonumber \\
&&S = \nu \overline{\left (1 + \lambda^2 \left ( \frac{\partial u}{\partial x} 
\right )^2 \right )^{1/2} \frac{\partial u}{\partial x}} 
= \nu \left [ K_4 + \lambda^2 K_5 \eta^2 \right ] \frac{\partial \bar{u}}{\partial x} 
+ O(\eta^4), 
\end{eqnarray}
and $K_4$, $K_5$ are defined by
\begin{eqnarray}
&&K_4 = \left (1 + \lambda^2 \left ( \frac{\partial \bar{u}}{\partial x} \right )^2 \right )^{1/2} \\
&&K_5 = \left (1 + \lambda^2 \left ( \frac{\partial \bar{u}}{\partial x} \right )^2 \right )^{-3/2} 
\Biggl [ \frac{1}{\gamma_T} \left ( \frac{\partial^2 \bar{u}}{\partial t \partial x} \right )^2 
\label{eq-K5} \\
&&\qquad {} + 
\frac{1}{\gamma_L} \left ( \frac{\partial^2 \bar{u}}{\partial x^2} \right )^2 \Biggr ]
\left (3 + 2\lambda^2 \left ( \frac{\partial \bar{u}}{\partial x} \right )^2 \right ) \nonumber
\end{eqnarray}

\subsection{Numerical comparison with some analytical solutions.}

Let us consider the following family of analytical solutions of the Burgers equation (\ref{eq-Burgers}):
\begin{equation} \label{sol-an-Burgers}
u = \frac{-2\nu}{\phi}\frac{\partial \phi}{\partial x}
\end{equation}
where 
\begin{eqnarray} \label{eq-phi}
&&\phi = A_0 + B_0 x + \left [ A_1\sin (\omega_1 x) + 
B_1\cos (\omega_1 x) \right ] \exp(-\nu \omega_1^2t) \\ 
&&\quad {} + \left [ A_2\sin (\omega_2 x) + B_2\cos (\omega_2 x) \right ] 
\exp(-\nu \omega_2^2t) \nonumber
\end{eqnarray}

Next, we shall use MacCormack scheme (see \cite{MacCormack}) to solve numerically problems (\ref{eq-Burgers}) and 
(\ref{eq-Burgers-prom})-(\ref{eq-K5}) and to compare the obtained results with the analytical solution 
(\ref{sol-an-Burgers}).

We shall compare with three different analytical solutions, obtained choosing the following three different sets of values for the constants of (\ref{sol-an-Burgers})-(\ref{eq-phi}) , $\eta$ and $\lambda$,
\begin{eqnarray}
&&\omega_1 = 2\pi, \omega_2 = 50\pi, A_0 = 1000, B_0 = -10,  
A_1=\frac{1}{\pi}, A_2 = \frac{1}{100\pi}, \label{set-ctes-1} \\
&&B_1 = B_2 = 0, \nu = \frac{1}{50000}, \eta = 0.1, \lambda =1.8e6, \nonumber 
\end{eqnarray}
\begin{eqnarray}
&&\omega_1 = 2\pi, \omega_2 = 50\pi, A_0 = 10, B_0 = -7,  
A_1=\frac{3}{\pi}, A_2 = \frac{3}{100\pi}, \label{set-ctes-2} \\
&&B_1 = B_2 = 0, \nu = \frac{1}{50000}, \eta = 0.1, \lambda =5000, \nonumber 
\end{eqnarray}
\begin{eqnarray}
&&\omega_1 = 2\pi, \omega_2 = 100\pi, A_0 = 1, B_0 = -10,  
A_1=\frac{-15}{\pi}, A_2 = \frac{1}{10\pi}, \label{set-ctes-3} \\
&&B_1 = \frac{-7}{\pi}, B_2 = \frac{-1}{100\pi}, \nu = \frac{1}{50000}, \eta = 0.01, \lambda =500, \nonumber 
\end{eqnarray}
where $\nu$ can be understood here as the viscosity, or as the inverse of Reynolds number in the non-dimensional version of (\ref{eq-Burgers}).

\begin{figure}[H]
  \begin{center}
  \vspace{-4.2cm}
  \makebox[\linewidth]{\includegraphics[width=.9\linewidth]{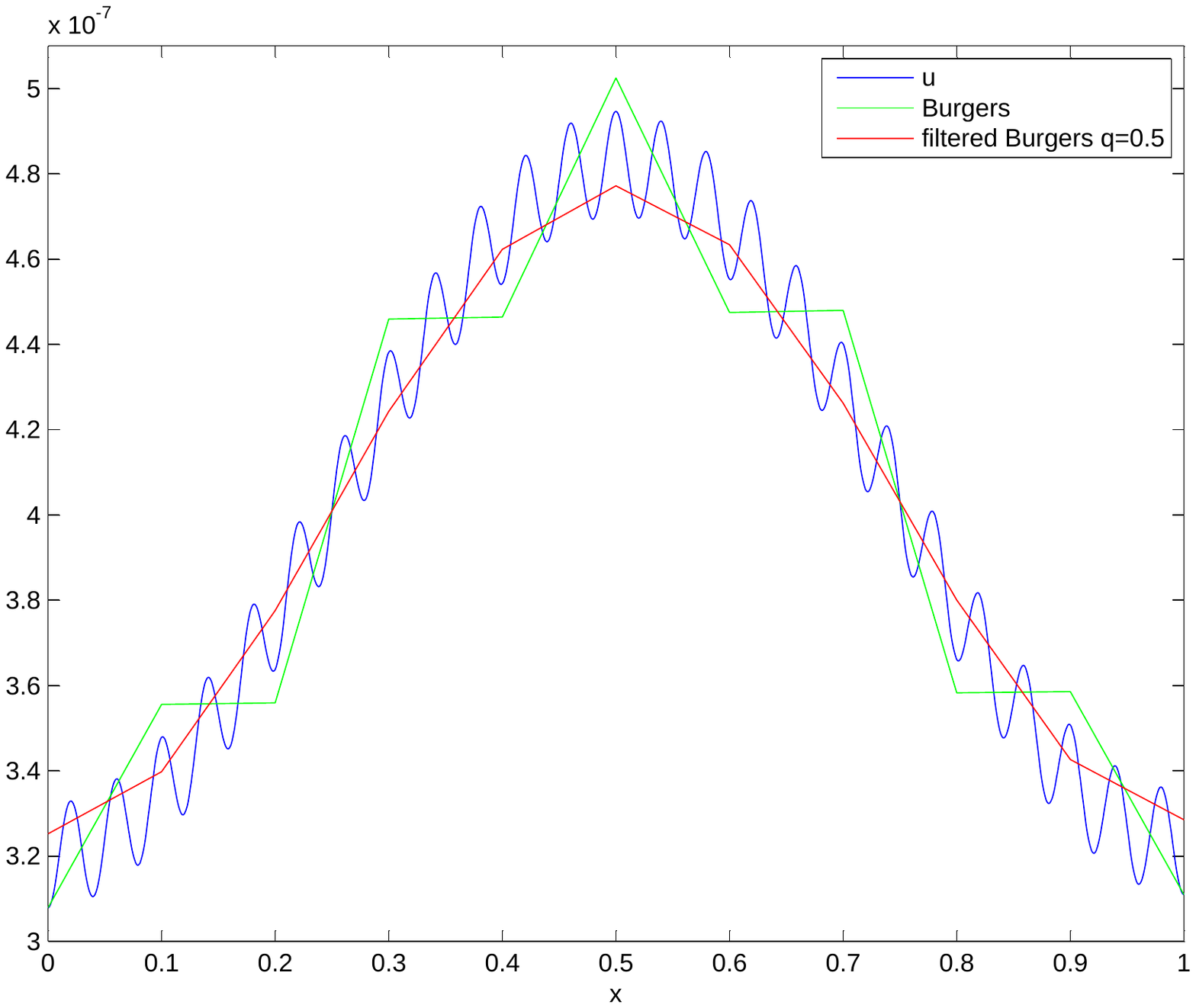}}
   \vspace{-4.2cm}
  \caption{Constants (\ref{set-ctes-1}), $\Delta x = 0.1$, $\Delta t = 0.01$, $t=1$.  \label{dibujo1}}
  \end{center}
\end{figure}

\begin{figure}[H]
  \begin{center}
  \vspace{-4.2cm}
  \makebox[\linewidth]{\includegraphics[width=.9\linewidth]{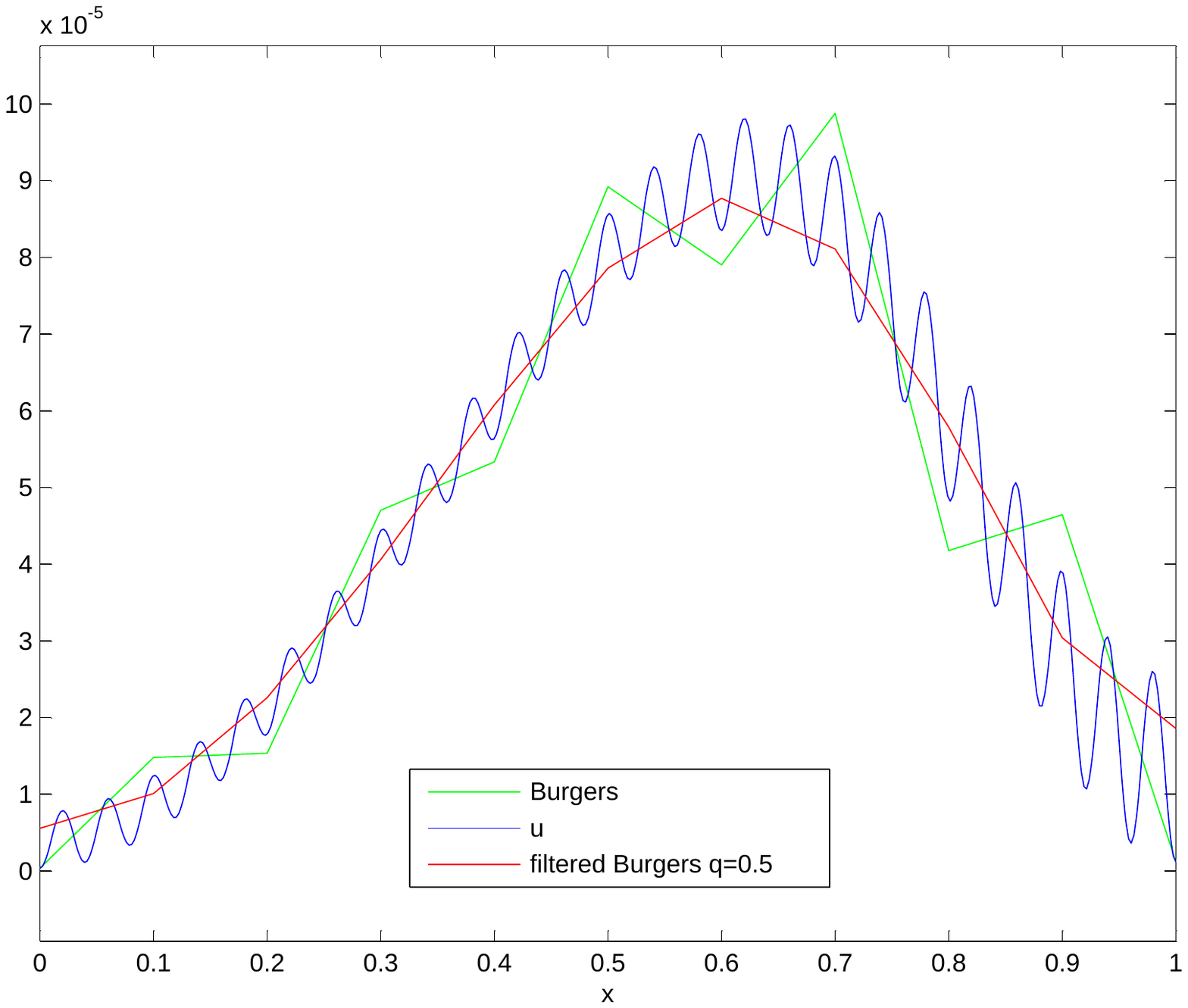}}
   \vspace{-4.2cm}
  \caption{Constants (\ref{set-ctes-2}), $\Delta x = 0.1$, $\Delta t = 0.01$, $t=1$. \label{dibujo2}}
  \end{center}
\end{figure}

\begin{figure}[H]
  \begin{center}
   \vspace{-4.2cm}
  \makebox[\linewidth]{\includegraphics[width=.9\linewidth]{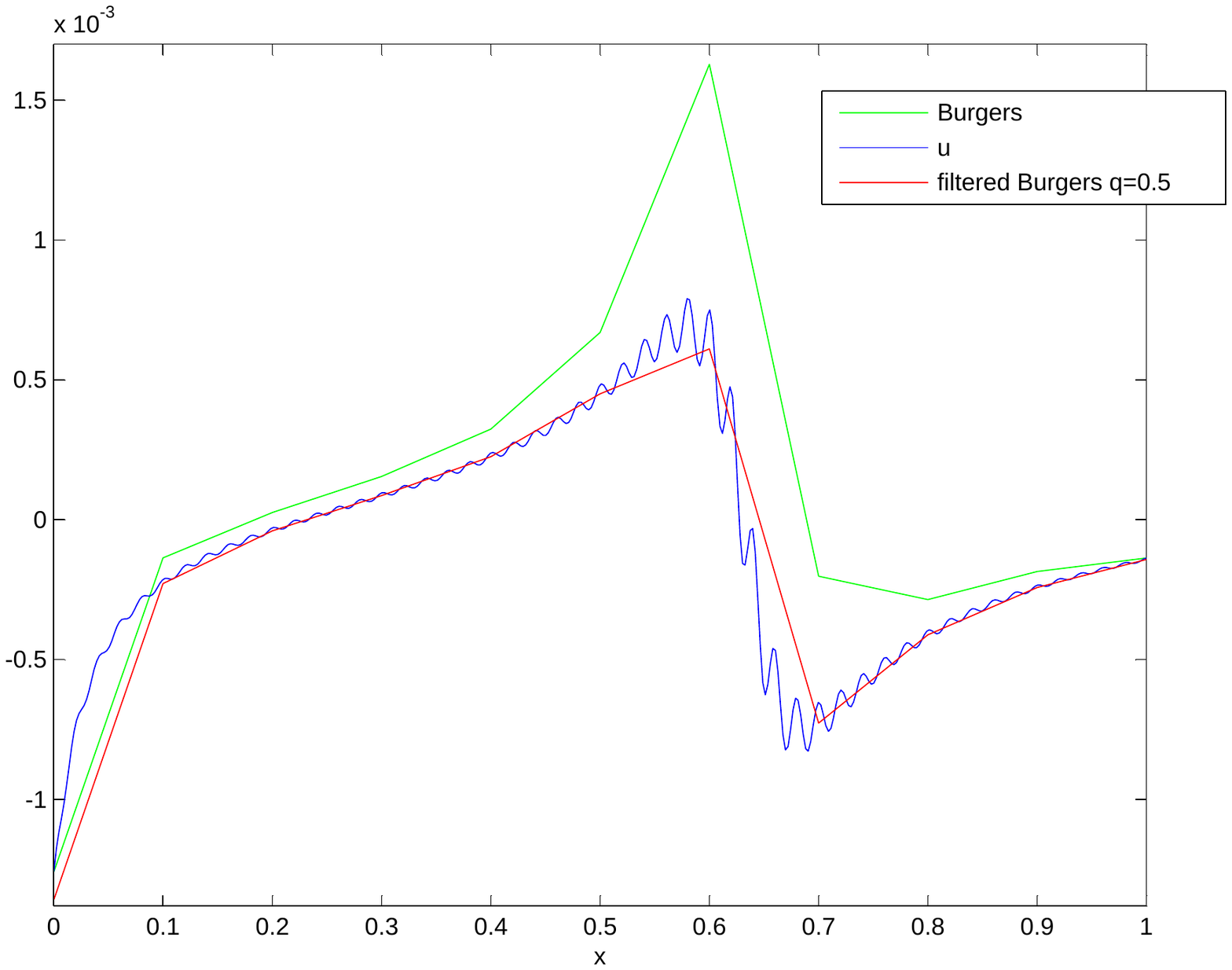}}
   \vspace{-4.2cm}
  \caption{Constants (\ref{set-ctes-3}), $\Delta x = 0.1$, $\Delta t = 0.01$, $t=1$. \label{dibujo3}}
  \end{center}
\end{figure}

We can see in figures \ref{dibujo1}-\ref{dibujo3} the good behavior of the numerical approximation of 
(\ref{eq-Burgers-prom})-(\ref{eq-K5}) if compared with the numerical approximation of (\ref{eq-Burgers}). For these simulations, we have chosen in (\ref{filtro-eq-Burgers}) $\gamma_T = \gamma_L = 6$. We have used a spatial step of $\Delta x = 0.1$ and a time step of $\Delta t = 0.01$. We have plotted in blue the exact solution, in green the numerical approximation of the Burgers equation (\ref{eq-Burgers}), and in red the numerical approximation of (\ref{eq-Burgers-prom})-(\ref{eq-K5}). We have represented the solutions at $t = 1$.

Despite we have taken a relatively large discretization step, model (\ref{eq-Burgers-prom})-(\ref{eq-K5}) gives in figures \ref{dibujo1}-\ref{dibujo3} a good approximation of the analytical solution of (\ref{eq-Burgers}), that is, we are able to approximate the ``large scale'' behavior of the solution without needing a small discretization step.

\begin{figure}[H]
  \begin{center}
   \vspace{-4.2cm}
  \makebox[\linewidth]{\includegraphics[width=.9\linewidth]{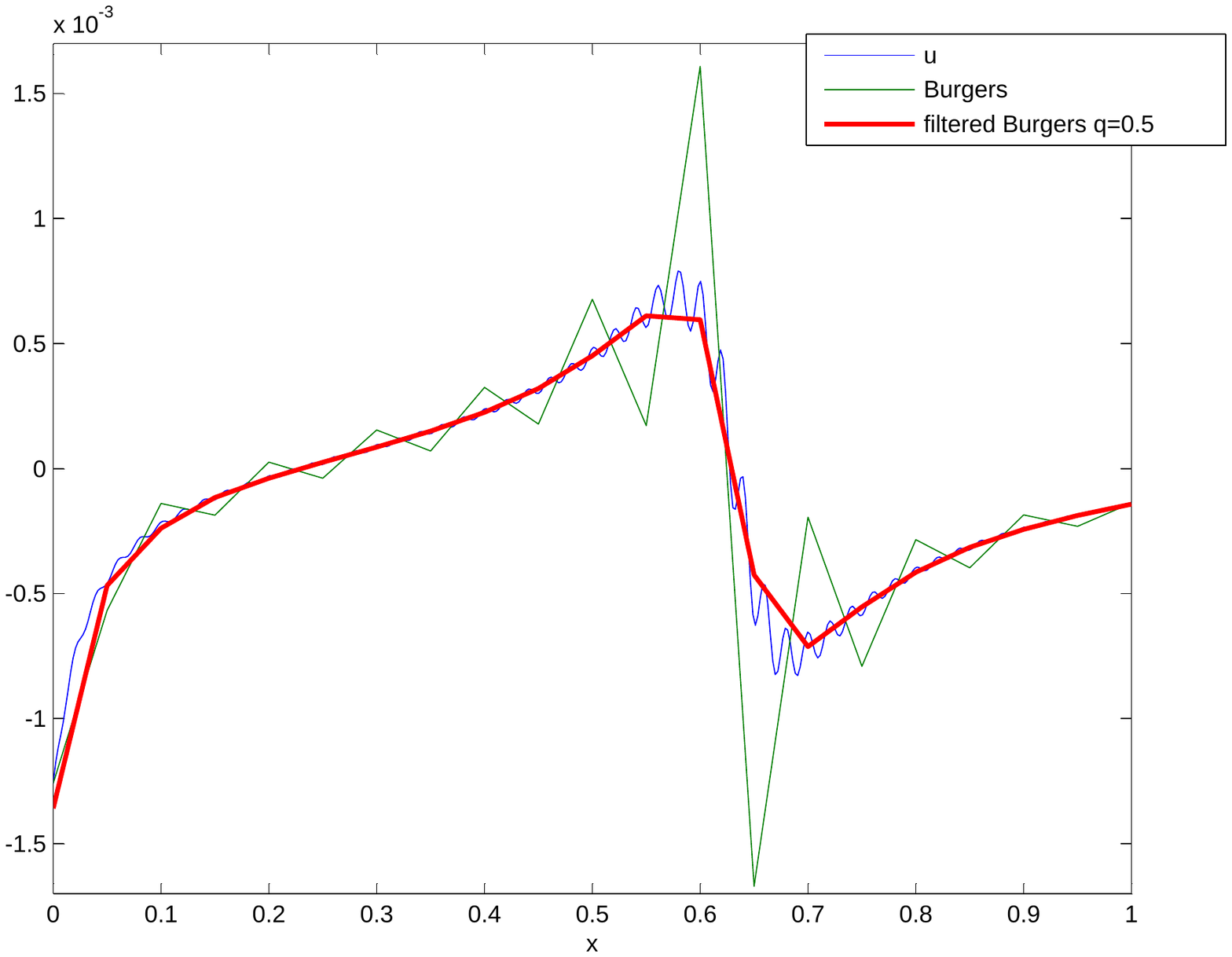}}
   \vspace{-4.2cm}
  \caption{Constants (\ref{set-ctes-3}), $\Delta x = 0.05$, $\Delta t = 0.005$, $t=1$. \label{dibujo4}}
  \end{center}
\end{figure}

\begin{figure}[H]
  \begin{center}
   \vspace{-4.2cm}
  \makebox[\linewidth]{\includegraphics[width=.9\linewidth]{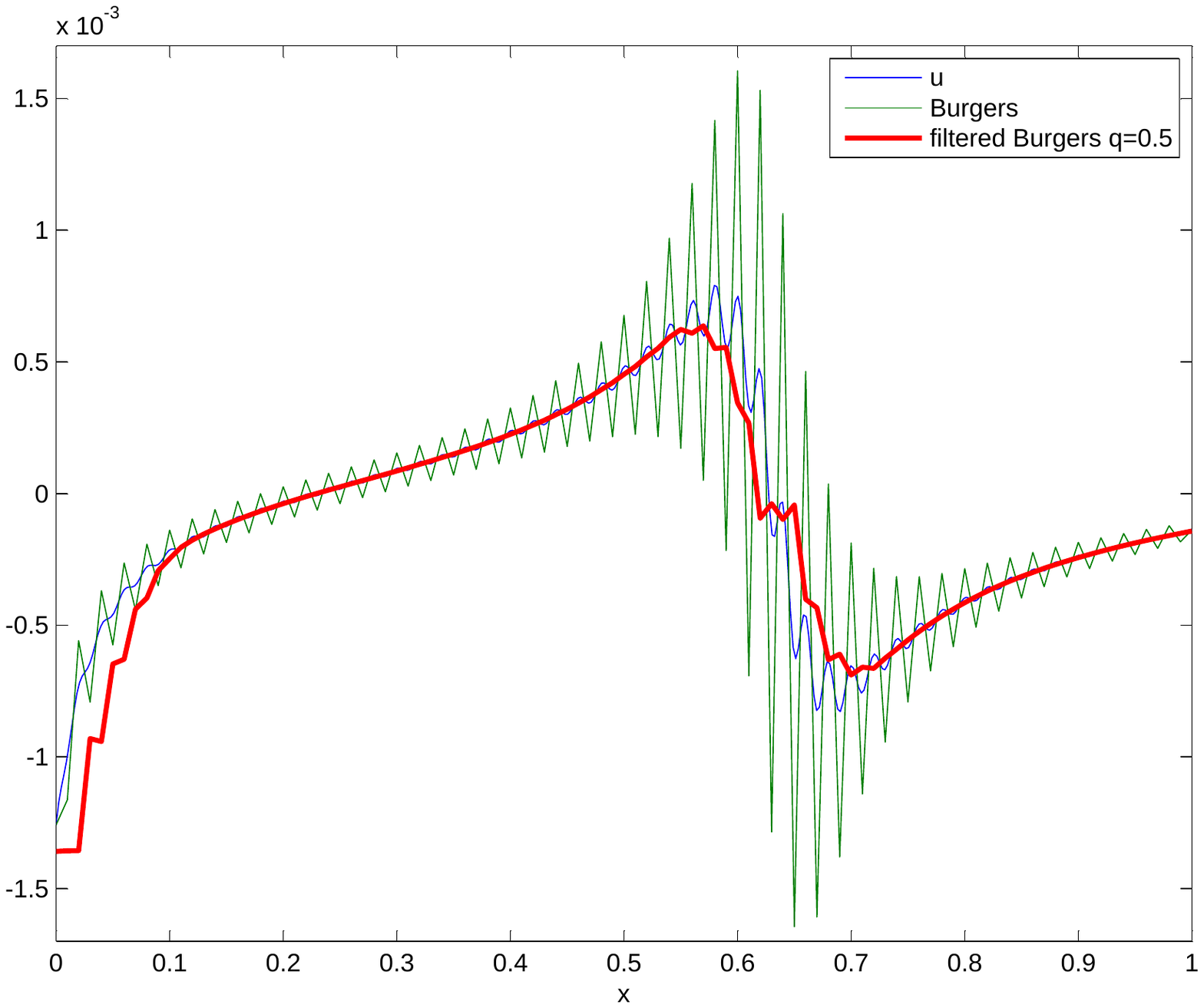}}
   \vspace{-4.2cm}
  \caption{Constants (\ref{set-ctes-3}), $\Delta x = 0.01$, $\Delta t = 0.001$, $t=1$. \label{dibujo5}}
  \end{center}
\end{figure}

Next, let us show what happens when we decrease the discretization step. In figures \ref{dibujo4} and \ref{dibujo5} we plot the analytical solution and their approximations for the same set of values of the constants as in figure \ref{dibujo3}, but with smaller discretization steps 
($\Delta x = 0.05$, $\Delta t = 0.005$ in figure \ref{dibujo4} and $\Delta x = 0.01$, $\Delta t = 0.001$ in figure \ref{dibujo5}). 

We can see that the MacCormack scheme has a bad numerical behavior when approximating Burgers equation (\ref{eq-Burgers}) (oscillations increase when decreasing discretization step), while its behavior when approximating equations (\ref{eq-Burgers-prom})-(\ref{eq-K5}) is very good.

\begin{figure}[H]
  \begin{center}
   \vspace{-4.2cm}
  \makebox[\linewidth]{\includegraphics[width=.9\linewidth]{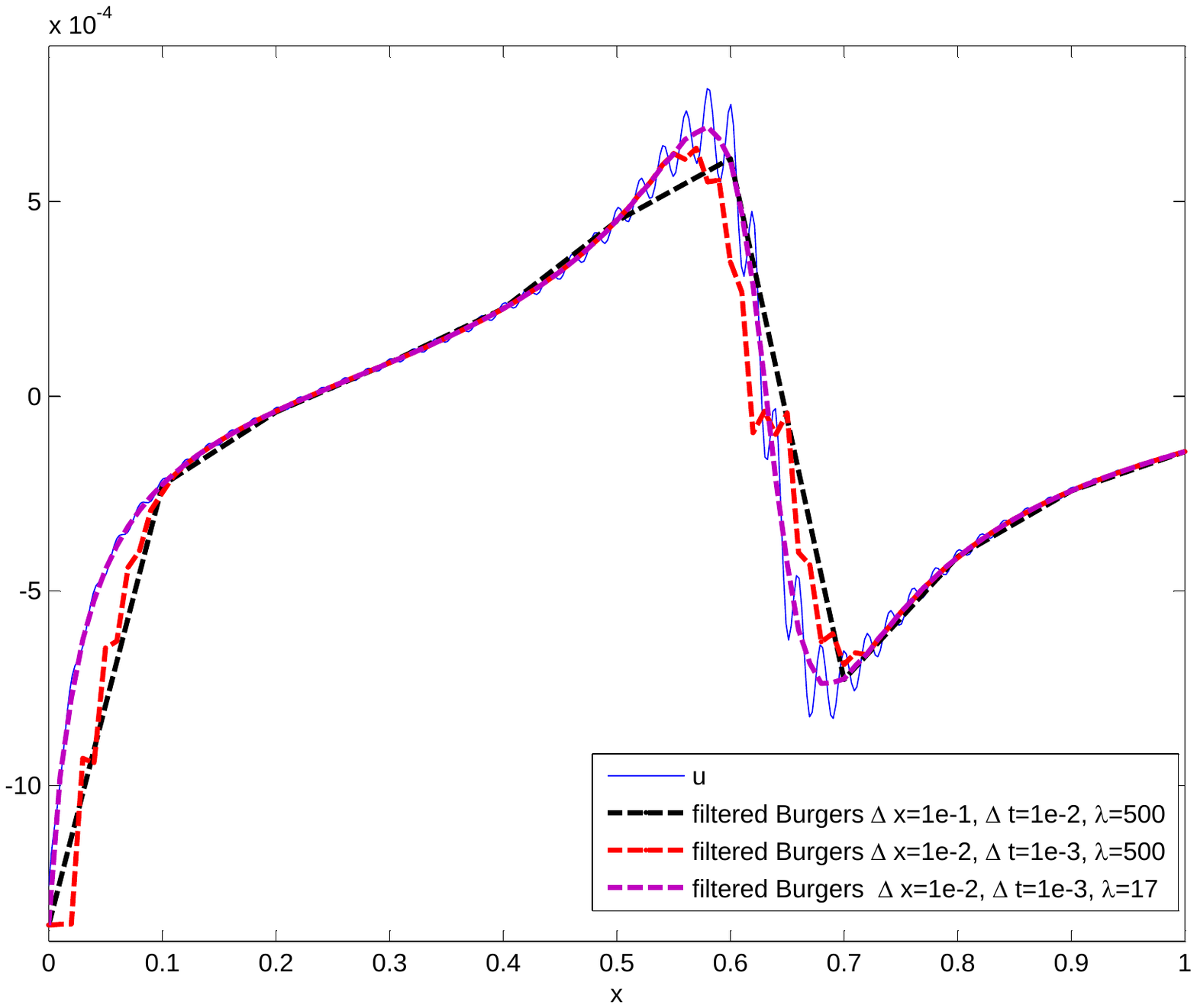}}
   \vspace{-4.2cm}
  \caption{Constants (\ref{set-ctes-3}), different values of $\lambda$. \label{dibujo6}}
  \end{center}
\end{figure}

Finally, we plot in figure \ref{dibujo6} the same analytical solution (set of values of the constants 
(\ref{set-ctes-3})) and several approximations of equations (\ref{eq-Burgers-prom})-(\ref{eq-K5}) with different values of $\lambda$ and discretization steps. This figure shows how the accuracy of the approximations depends not only on the discretization step, but also on the election of $\lambda$.

\subsection{Another numerical example}

\begin{figure}[H]
  \begin{center}
   \vspace{-1.5cm}
  \makebox[\linewidth]{\includegraphics[width=.86\linewidth]{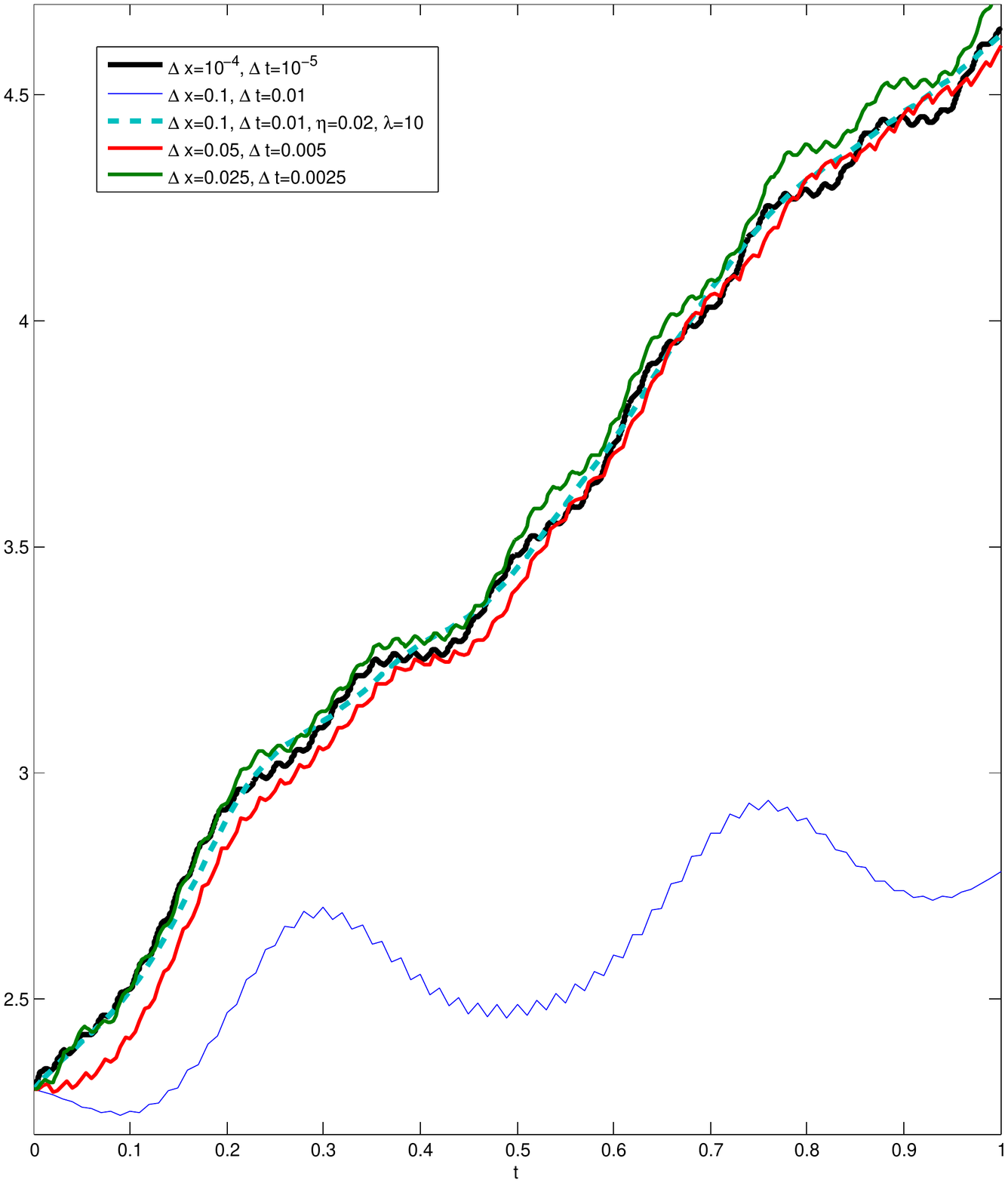}}
   \vspace{-2.5cm}
  \caption{Different approximations at $x=0.5$. \label{dibujo7}}
  \end{center}
\end{figure}

\begin{figure}[H]
  \begin{center}
   \vspace{-1.5cm}
  \makebox[\linewidth]{\includegraphics[width=.9\linewidth]{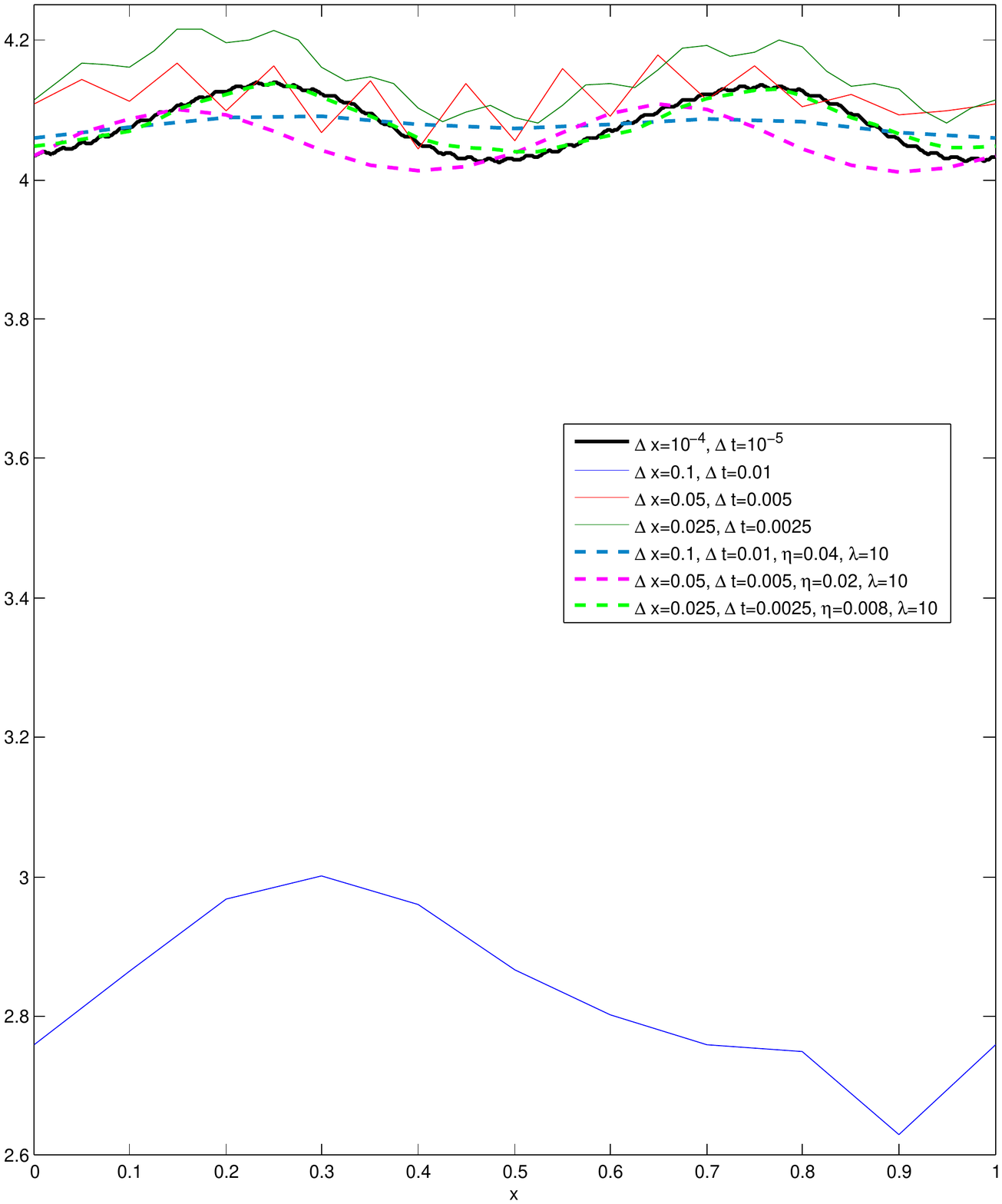}}
   \vspace{-2.5cm}
  \caption{Different approximations at $t=0.7$. \label{dibujo8}}
  \end{center}
\end{figure}

Now, we shall use again the MacCormack scheme 
to solve numerically equation (\ref{eq-Burgers}) with periodic boundary conditions at 
$x=0$ and $x=1$, a constant initial velocity and a periodic external force. We show in figures 
\ref{dibujo7}-\ref{dibujo8} the results obtained taking $u(0,x) = 2.3$, 
 $\nu = 1/5000$ and the external force $f$ given by
\begin{equation}
f = 2.3 + \sin (4\pi x) + 2.9\sin (99\pi x) + \sin (4\pi t) + 2.9\sin (99\pi t)
\end{equation}

In figure \ref{dibujo7} we have plotted different approximations of the solution at point $x=0.5$. 
We have used four different different discretization steps ($\Delta x = 0.1$, $\Delta x = 0.05$, 
$\Delta x = 0.025$, $\Delta x = 10^{-4}$). In order to compare the accuracy of the filtered model, 
let us consider that the ``exact'' solution is the one given by the smallest discretization step (we have plotted this solution in black in figures \ref{dibujo7}-\ref{dibujo8}). We can see then that the filtered solution (computed with the larger step size, $\Delta x = 0.1$, and 
plotted with a discontinuous line) is more accurate than the approximations obtained by solving 
Burgers equation (\ref{eq-Burgers}) with smaller step sizes.

We have plotted in figure \ref{dibujo8} some approximations of the solution, at the instant $t=0.7$, 
for the discretization steps mentioned above ($\Delta x = 0.1$, $\Delta x = 0.05$, 
$\Delta x = 0.025$). We can observe that filtered approximations (discontinuous lines) are closer to the ``exact'' solution ($\Delta x = 10^{-4}$) than the non filtered approximations 
(continuous lines) with the same, and even smaller, discretization step. 
This is specially clear for the coarser 
discretization step ($\Delta x = 0.1$), where the filtered approximation 
essentially captures the exact solution, while the unfiltered approximation remains far from it.

\section{Conclusions}

A new LES model has been deduced. We have applied a filter to the generalized Navier-Stokes equations with a nonlinear effective viscosity and then we have approximated the subgrid-scale tensor using the Clark approximation. The new model thus obtained reminds us the dynamic procedure of Germano.

We have performed some numerical simulations to test this new model. We have solved Burgers equation and the model (\ref{eq-Burgers-prom})-(\ref{eq-K5}), obtained with the same technique used previously for Navier-Stokes equations. We have compared the results with analytical and numerical solutions of the Burgers equations, and we have seen the good numerical behavior of the model that we have deduced. Even for large discretization steps the filtered model provides quite good approximations of the exact solutions.

\section*{Acknowledgements}
This research has been partially supported by Ministerio de Economía y Competitividad under grant MTM2012-36452-C02-01 with the participation of FEDER.


\begin{thebibliography}{99}
\footnotesize

\bibitem{Argyropoulos} \textrm{C.\ D.\ Argyropoulos, N.\ C.\ Markatos}, \titleref{Recent advances on the numerical modelling of turbulent flows.} Applied Mathematical Modelling \textbf{39} (2015), 693--732.

\bibitem{Layton-libro} \textrm{L.\ C.\ Berselli, T.\ Iliescu, W.\ J.\ Layton}, \titleref{Mathematics of Large Eddy Simulation of Turbulent Flows.} Springer, 2006.

\bibitem{Carati1} \textrm{D.\ Carati, F.\ S.\ Winckelmans, H.\ Jeanmart}, \titleref{On the modelling of the subgrid-scale and filtered-scale stress tensors in large-eddy simulation.} J. Fluid Mech. \textbf{441} (2001), 119--138.

\bibitem{Carati2} \textrm{D.\ Carati, A. A. Wray}, \titleref{Time filtering in large eddy simulation.} Center for Turbulence Research, Proceedings of the Summer Program 2000.

\bibitem{Dakhoul1} \textrm{Y.\ M.\ Dakhoul, K.\ W.\ Bedford}, \titleref{Improved averaging method for turbulent flow simulation. Part I: Theoretical development and application to Burgers' transport equation.} International Journal for Numerical Methods in Fluids \textbf{6} (1986), 49--64.

\bibitem{Dakhoul2} \textrm{Y.\ M.\ Dakhoul, K.\ W.\ Bedford}, \titleref{Improved averaging method for turbulent flow simulation. Part II: Calculations and verification.} International Journal for Numerical Methods in Fluids \textbf{6} (1986), 65--82.

\bibitem{Germano} \textrm{M.\ Germano, U.\ Piomelli, P.\ Moin, W.\ H.\ Cabot}, \titleref{A dynamic subgrid-scale eddy viscosity model.} Phys. Fluids A \textbf{3} (1991), 1760--1765.

\bibitem{Guermond2} \textrm{J.-L.\ Guermond, S.\ Prudhomme}, \titleref{Mathematical analysis of a spectral hyperviscosity LES model for the simulation of turbulent flows.} Mathematical Modelling and Numerical Analysis \textbf{37}, n. 6, (2003), 893--908.

\bibitem{Guermond} \textrm{J.-L.\ Guermond, J.\ T.\ Oden, S.\ Prudhomme}, \titleref{Mathematical Perspectives on Large Eddy Simulation Models.} J. Math. Fluid Mech. \textbf{6} (2004), 194--248.

\bibitem{Lesieur} \textrm{M.\ Lesieur, O.\ M\'etais}, \titleref{New trends in large-eddy simulations of turbulence.} Annu. Rev. Fluid Mech. \textbf{28} (1996), 45--82.

\bibitem{Lilly} \textrm{D.\ K.\ Lilly}, \titleref{A proposed modification of the Germano subgrid-scale closure method.} Phys. Fluids A \textbf{4} (1992), 633--635.

\bibitem{MacCormack} \textrm{R.\ W.\ MacCormack}, \titleref{Numerical solution of the interaction of a shock wave with a laminar boundary layer.} In: M. Holt (Ed.), Proceedings 2nd Int. Conf. on Num. Methods in Fluid Dynamics, Springer-Verlag, 1971, pp. 151--163.

\bibitem{Pruett} \textrm{C.\ Pruett}, \titleref{Temporal large-eddy simulation: theory and implementation.} Theor. Comput. Fluid Dyn. \textbf{22} (2008), 275--304.

\bibitem{Sagaut} \textrm{P.\ Sagaut}, \titleref{Large Eddy Simulation for Incompressible Flows.} Springer, 2006.

\bibitem{Vreman} \textrm{B.\ Vreman, B.\ Geurts, H.\ Kuerten}, \titleref{Large-Eddy Simulation of the Temporal Mixing Layer Using the Clark Model.} Theoret. Comput. Fluid Dynamics \textbf{8} (1996), 309--324.

\end{thebibliography}
\end{document}